\documentclass[a4,12pt]{article}
\usepackage{amsmath}
\usepackage{amsfonts}
\usepackage{amssymb}

\def\Diff{\mathrm{Diff}}
\def\Vect{\mathrm{Vect}}
\newtheorem{thm}{Theorem}

\title{Integrable invariant Sobolev metrics on the Abelian extension
of the diffeomorphism group of the circle and two-component generalizations 
of the Camassa-Holm equation}
\author{P.~A.~Kuzmin\\
{\em Department of Physics, Moscow State University,} \\ {\em
Leninskie gory, Moscow 119992, Russia}\\ {\it e-mail:
kuzminp@list.ru} }

\date{ \ \ \ }
\begin{document}
\maketitle
\begin{abstract}
In this note we classify some
integrable invariant Sobolev metrics on the Abelian extension
of the diffeomorphism group of the circle.
We also derive a new two-component generalization
of the Camassa-Holm equation. The system obtained appears to
be unique bi-Hamiltonian flow on the coadiont obrit of
$\Diff_+(S^1)\ltimes C^\infty(S^1)$ generalizes the Camassa-Holm
flow on $\Vect(S^1)^*$.
\end{abstract}

It is known that the Camassa-Holm equation can be interpreted
as the flow on the coadjoint orbit of the diffeomorphism
group of the circle. The Camassa-Holm equation was derived in \cite{CH,CHH,FF}.
Two-component generalizations of the  equation were  introduced in
\cite{LZ,OlRo,AGZ,Fa,CLZ}, where some properties of the systems
were studied. All these systems appear to be bi-Hamiltonian
flows on the regular part of the dual space of the Lie algebra
$\frak g=\Vect(S^1)\ltimes C^\infty(S^1)$.

In this paper we consider a class of the Hamitlonian flows on the
space $\frak g^*_{reg}$ with quadratic Hamiltonian functions of special form.
We are going to classify all the flows those are bi-Hamiltonian
relative to a modified Lie-Poisson structure on $\frak g^*_{reg}$.
We derive a family of two-component completely integrable  generalizations
of the Camassa-Holm equation.

There is a standard way to construct a pair of compatible Poisson
structures and obtain an integrable flow.
Let $\frak g$ be a Lie algebra, $\frak g^*$ be the dual
space of $\frak g$ (the regular part of the dual space of $\frak g$).
One can define the canonical Lie-Poisson structure on $\frak g^*$
given by
\begin{equation}\label{twoast}
\{f,g\}(m)=m([df,dg])
\end{equation}
for all smooth $f, g: \frak g^*\to\mathbb R$.
Each cocycle $\omega\in Z^2(\frak g)$ defines another Poisson
structure on $\frak g^*$, which is compatible with the canonical
Lie-Poisson structure. This structure is called the modified Lie-Poisson
structure.

Let us recall some facts on the Camassa-Holm equation.
The canonical Lie-Poisson structure on the regular part of the
dual space of $\Vect(S^1)$ is represented by
the skew-symmetric operator
$$
j_0(m)=m_x+2m\partial,\,\,m\in\Vect(S^1)^*_{reg}.
$$
The regular $\frak g^*_{reg}$ part of the dual space is the set of all
linear functionals on $\frak g=\Vect(S^1)$ of the form
$$
\bar m:\,\,f\mapsto\int_{S^1}fm\,dx,\,\,\,m\in C^\infty(S^1),
$$
and we identify the functional $\bar m$ with the function $m$. \cite{Ki1}
For the modified Lie-Poisson structure,
which is obtained from the  Virasoro cocycle, one has
$$
j(m)=2m_0\partial+(m_0)_x-c_1\partial^3,
$$
where $m_0\in\Vect(S^1)^*$, $c_1\in\mathbb R$, and $\partial=\frac{d}{dx}$.
Consider the Hamiltonian vector field $x(m)=j_0(m)\delta h_0(m)$ on $\Vect(S^1)^*_{reg}$
with  Hamiltonian function
$$
h_0(m)=\frac{1}{2}\int\limits_{S^1}m(1-\partial^2)^{-1}m\,dx.
$$
Here $\delta \phi(m)$ denotes the gradient of the smooth functional $\phi(m)$.
The field $x(m)$ is known to be Hamiltonian with respect to the Poisson
structure $j_1(m)=\partial-\partial^3$ with Hamiltonian function
$$
h_1(m)=\frac{1}{2}\int\limits_{S^1}(u^3+uu_x^2)dx,\,\, u=(1-\partial^2)^{-1}m,
$$
i.e. $x(m)=j_0(m)\delta h_0(m)=j_1(m)\delta h_1(m)$.
The corresponding evolution equation $m_t+x(m)=0$ is called the
Camassa-Holm equation:
$$
m_t+m_xu+2mu_x=0,\,\,u=\left(1-\partial^2\right)^{-1}m.
$$

The main goal of this note is to construct  analogous bi-Hamiltonian
systems on the space $\frak g^*_{reg}$, where $\frak g=\Vect(S^1)\ltimes C^\infty(S^1)$,
and $C^\infty(S^1)$ is the space of all (smooth) tensor densities on $S^1$ of degree $0$,
i.e. the $\Vect(S^1)$-module structure on $C^\infty(S^1)$ is defined as
$$
L_{f}:\,\,C^\infty(S^1)\to C^\infty(S^1),\,\, a(x)\mapsto f(x)a^\prime(x),
$$
where $f\in\Vect(S^1)$.
The canonical Lie-Poisson structure on $\frak g^*_{reg}$
is given by the operator
\begin{equation}\label{thrast}
J_0\left(\left(\begin{matrix} m \\ p \end{matrix}\right)\right)=
\left(\begin{matrix} m_x+2m\partial & p\partial \\
p_x+p\partial & 0 \end{matrix}\right),\,\,
\left(\begin{matrix} m \\ p \end{matrix}\right)\in\frak g^*_{reg},
\end{equation}
which defines the coadjoint action of the semi-direct product
$\Vect(S^1)\ltimes C^\infty(S^1)$. \cite{MOR}
It is known that $H^2(\Vect(S^1)\ltimes C^\infty (S^1))=\mathbb R^3$
(see, for example, \cite{OR1,ACKP}), so the modified Lie-Poisson structure
is represented by the constant  (i.e. independent on the point
$\left(\begin{matrix} m \\ p \end{matrix}\right)\in\frak g_{reg}^*$)
operator
\begin{equation}\label{mps}
J=\left(\begin{matrix}(m_0)_x+2m_0\partial-c_1\partial^3 & p_0\partial+c_2\partial^2 \\
(p_0)_x+p_0\partial-c_2\partial^2 & 2c_3\partial \end{matrix}\right),
\end{equation}
where $\left(\begin{matrix} m_0 \\ p_0 \end{matrix}\right)\in\frak g^*_{reg}$ and
$c_1,\,c_2,\,c_3$ are constants.
Note that the structures (\ref{thrast}), (\ref{twoast}) are actually
defined on the coadjoint orbit
$(\Diff_+(S^1)\ltimes C^{\infty}(S^1))/(S^1\times S^1)$ containing a constant
moment $\left(\begin{matrix}m_0\\p_0\end{matrix}\right)\in\frak g^*_{reg}$,
$m_0(x)=const,\,p_0(x)=const$ \cite{Ki1}.

Consider a class of quadratic Hamiltonian functions on $\frak g^*_{reg}$ given
by
\begin{equation}\label{ham}
H_0\left(\left(\begin{matrix} m \\ p \end{matrix}\right)\right)=\frac{1}{2}
\left\langle \left(\begin{matrix} m \\ p \end{matrix}\right),M^{-1}
\left(\begin{matrix} m \\ p \end{matrix}\right)\right\rangle,
\end{equation}
where $\langle,\rangle$ denotes the $L^2$ pairing
$$
\left\langle \left(\begin{matrix} m_1 \\ p_1 \end{matrix}\right),
\left(\begin{matrix} m_2 \\ p_2 \end{matrix}\right)\right\rangle=
\int\limits_{S^1}\left(m_1(x)m_2(x)+p_1(x)p_2(x)\right)dx,
$$
and the symmetric "inertia" operator is chosen to be of the form
\begin{equation}\label{inert}
M=\left(\begin{matrix} \varepsilon\left(1-\partial^2\right)+\sum\limits_{i=2}^{n_0} a_i\partial^{2i} &
 \sum\limits_{i=0}^{n_1} (-1)^i b_i\partial^{i}\\
 \sum\limits_{i=0}^{n_1} b_i\partial^{i} & \gamma \end{matrix}\right),
\end{equation}
where $a_i,\,b_i,\,\gamma,\,\varepsilon$ are  constants such that $M$ is
strictly positive and invertible.
This Hamiltonian function
can be considered as a natural generalization of the Hamiltonian function $h_0(m)$.
Now we state the main theorem of this note:
\begin{thm}
Assume that $m_0,\,p_0$ (see (\ref{mps})) are constant functions.
The only vector field of the form $X=J_0\delta H_0$
(where $J_0,\,H_0$ are given by (\ref{thrast}), (\ref{ham})) that is
bi-Hamiltonian relative to the modified Lie-Poisson structure (\ref{mps})
is the Hamiltonian vector field defined by the following "inertia"
operator:
\begin{equation}\label{ast}
M= \left(\begin{matrix} \varepsilon\left(1-\partial^2\right) & \alpha-\beta\partial \\
\alpha+\beta\partial & \gamma \end{matrix}\right),
\end{equation}
where $\alpha,\,\beta,\,\gamma,\,\varepsilon$ are constants.
\end{thm}

We  obtain the second Poisson structure
\begin{equation}\label{sps}
J_1= \left(\begin{matrix} \varepsilon\left(\partial-\partial^3\right) & \alpha\partial-\beta\partial^2 \\
\alpha\partial+\beta\partial^2 & \gamma\partial \end{matrix}\right),
\end{equation}
and the second Hamiltonian function
$$
H_1\left(\left(\begin{matrix} m \\ p \end{matrix}\right)\right)=
\frac{1}{2}\int\limits_{S^1}(2\alpha u^2v+2\beta uu_xv+\gamma uv^2+\varepsilon u^3+\varepsilon uu_x^2)dx,
$$
$$
\left(\begin{matrix} u \\ v \end{matrix}\right)=M^{-1}\left(\begin{matrix} m \\ p \end{matrix}\right),
$$
such that $X=J_0\delta H_0=J_1\delta H_1$.
The equations obtained have the form
$$
\left\{\begin{array}{lcl}
m_t+m_xu+2mu_x+pv_x=0\\
p_t+(pu)_x=0,\end{array}\right.
$$
where $m=\varepsilon u-\varepsilon u_{xx}+\alpha v-\beta v_x,\,\,
p=\alpha u+\beta u_x+\gamma v$.

The proof of the theorem is based on the following simple observation.
This method was developed in \cite{CK,CKL}.
Let $X$  be a smooth vector field on $\frak g^*_{reg}$. Suppose that
$X$ is Hamiltonian with respect to the constant  (i.e. independent on the point
of $\frak g_{reg}^*$)  Poisson structure $J$.
This means that there is exists a smooth function $H:\frak g^*_{reg}\to\mathbb R$
such that
$$
X=J\delta H(m).
$$
The operator $d\delta H(m)$ is symmetric, since $H(m)$ is a smooth function.
Therefore, the operator
$$
S(m)=J(d\delta H(m))J=d(J\delta H(m))J=(dX)J
$$
is symmetric.
We obtained a necessary condition for $X$ to be  Hamiltonian with
respect to the Poisson structure $J$.

We restrict our considerations to the  modified Lie-Poisson
structures (\ref{mps}) such that $m_0,\,p_0$ are constants.
Thus we obtain the following condition
for $X$ to be Hamiltonian with respect to $J$:
\begin{equation}\label{cond}
\text{The operator}\,\,
(dX)J=
\left(\begin{matrix} u\partial & v_x\\
0 & u\partial+u_x
\end{matrix}\right)J+
\left(\begin{matrix} m_x+2m\partial & p\partial \\
p\partial+p_x & 0
\end{matrix}\right)M^{-1}J\,\, \text{is symmetric,}
\end{equation}
where
$$
J=\left(\begin{matrix}2m_0\partial-c_1\partial^3 & p_0\partial+c_2\partial^2 \\
p_0\partial-c_2\partial^2 & 2c_3\partial \end{matrix}\right),\,\,
\left(\begin{matrix} u \\ v \end{matrix}\right)
=M^{-1}\left(\begin{matrix} m \\ p \end{matrix}\right).
$$
Note that the (operator) expression
\begin{equation}
((dX)J)\left(\left(\begin{matrix}m\\p\end{matrix}\right)\right)=
\left(\begin{matrix} u\partial & v_x\\
0 & u\partial+u_x
\end{matrix}\right)J+
\left(\begin{matrix} m_x+2m\partial & p\partial \\
p\partial+p_x & 0
\end{matrix}\right)M^{-1}J
\end{equation}
is linear in  $\left(\begin{matrix}m\\p\end{matrix}\right)$.
We are now in a position to obtain a polynomial representation for
the following trilinear functionals:
\begin{equation}\label{matel}
\left\langle\left(\begin{matrix}m_1\\p_1\end{matrix}\right),
((dX)J)\left(\left(\begin{matrix}m\\p\end{matrix}\right)\right)
\left(\begin{matrix}m_2\\p_2\end{matrix}\right)\right\rangle
\,\,\text{and}\,\,
\left\langle((dX)J)\left(\left(\begin{matrix}m\\p\end{matrix}\right)\right)
\left(\begin{matrix}m_1\\p_1\end{matrix}\right),
\left(\begin{matrix}m_2\\p_2\end{matrix}\right)\right\rangle.
\end{equation}
We can achieve this substituting
$$
\left(\begin{matrix}m\\p\end{matrix}\right)=
\left(\begin{matrix}e^{ikx}\\e^{ilx}\end{matrix}\right),\,\,\,
\left(\begin{matrix}m_1\\p_1\end{matrix}\right)=
\left(\begin{matrix}e^{ik_1x}\\e^{il_1x}\end{matrix}\right),\,\,\,
\left(\begin{matrix}m_2\\p_2\end{matrix}\right)=
\left(\begin{matrix}e^{ik_2x}\\e^{il_2x}\end{matrix}\right).
$$
One can see that expressions (\ref{matel}) vanish unless
$l=k,\,l_1=k_1,\,l_2=k_2,\,k=-k_1-k_2$.
Therefore
\begin{eqnarray}\label{pc}
\left\langle\left(\begin{matrix}e^{ik_1x}\\e^{ik_1x}\end{matrix}\right),
((dX)J)\left(\left(\begin{matrix}e^{-ik_1x-ik_2x}\\e^{-ik_1x-ik_2x}\end{matrix}\right)\right)
\left(\begin{matrix}e^{ik_2x}\\e^{ik_2x}\end{matrix}\right)\right\rangle, \nonumber \\
\left\langle((dX)J)\left(\left(\begin{matrix}e^{-ik_1x-ik_2x}\\e^{-ik_1x-ik_2x}\end{matrix}\right)\right)
\left(\begin{matrix}e^{ik_1x}\\e^{ik_1x}\end{matrix}\right),
\left(\begin{matrix}e^{ik_2x}\\e^{ik_2x}\end{matrix}\right)\right\rangle
\end{eqnarray}
are polynomials in $k_1$ and $k_2$.
Straightforward calculation yields the result posed:
the polonomials (\ref{pc}) coincide and (\ref{cond}) is satisfied
only if the "inertia" operator (\ref{inert}) has the form (\ref{ast}).
The second Poisson structure turns to be of the form (\ref{sps}).

Let us consider some special cases of the bi-Hamiltonian flow constructed.
Setting $\alpha=\beta=0,\,\varepsilon=-\gamma=1$ we recover the system
studied in \cite{CLZ,AGZ}:
\begin{equation}\label{ch21}
\left\{\begin{array}{lcl}
m_t+m_xu+2mu_x-pp_x=0\\
p_t+(pu)_x=0,\end{array}\right.
\end{equation}
where $m=u-u_{xx}$.
If, in turn, we set $\varepsilon=\gamma=0,\,\alpha=\beta=1$, we
recover the system introduced in \cite{LZ,Fa}:
\begin{equation}\label{ch22}
\left\{\begin{array}{lcl}
(v-v_x)_t+(2uv-v_xu)_x=0\\
(u+u_x)_t+(u^2+uu_x)_x=0,\end{array}\right.
\end{equation}
which is equivalent to the system considered in \cite{Fa} under
$x\mapsto -x,\,u\mapsto v,\,v\mapsto u$.
The same bi-Hamiltonian structures for (\ref{ch21}),\,(\ref{ch22}) were
obtained in \cite{CLZ,AGZ,LZ,Fa}.

\end{document}